%% file: bpde-paper.tex
\documentclass[10pt]{article}

\usepackage[a4paper]{geometry}
\usepackage{indentfirst}
\usepackage[pdfpagemode=UseNone]{hyperref}
\usepackage{amssymb, amsmath, amsthm}

\newcommand{\mynewtheorem}[2]{\newtheorem{#1}{\indent #2}}
\newcommand{\myalttheorem}[2]{\newtheorem*{#1}{\indent #2}}
\mynewtheorem{question}{Question}
\myalttheorem{question*}{Question}
\mynewtheorem{conjecture}{Conjecture}
\myalttheorem{conjecture*}{Conjecture}
\mynewtheorem{lemma}{Lemma}
\mynewtheorem{corollary}{Corollary}
\myalttheorem{corollary*}{Corollary}
\mynewtheorem{proposition}{Proposition}
\myalttheorem{proposition*}{Proposition}
\mynewtheorem{theorem}{Theorem}
\myalttheorem{theorem*}{Theorem}
\newenvironment{myproof}[1][Proof]{\begin{proof}[\indent #1]}{\end{proof}}

\newcommand{\lhull}{[\![}
\newcommand{\rhull}{]\!]}
\newcommand{\vd}{\vol_\textnormal{Disc}}

\DeclareMathOperator{\vol}{vol}
\DeclareMathOperator{\mix}{mix}
\DeclareMathOperator{\sol}{sol}
\DeclareMathOperator{\add}{add}
\DeclareMathOperator{\disc}{disc}
\DeclareMathOperator{\side}{side}
\DeclareMathOperator{\sd}{SD}

\begin{document}

\title{\textbf{Bootstrap Percolation and\\ Partial Difference Equations}}
\author{Nikolai Beluhov}
\date{}

\maketitle

\begin{abstract} We study a bootstrap percolation process on $\mathbb{Z}^d$ where each newly occupied point completes a copy of one of the patterns in some fixed collection $\mathcal{A}$. We aim, given $\mathcal{A}$, to find the smallest size of a percolating seed. We show that, for two patterns in two dimensions, in essence the answer is given by the mixed volume of the patterns' convex hulls. Our motivation for studying this setup comes from the theory of partial difference equations. Conversely, our considerations yield one consistency criterion for systems of two partial difference equations in two dimensions. \end{abstract}

\input{bpde-01-intro.tex}

\input{bpde-02-prelim.tex}

\input{bpde-03-plane.tex}

\input{bpde-04-cons-i.tex}

\input{bpde-05-cons-ii.tex}

\input{bpde-06-sd.tex}

\input{bpde-07-further.tex}

\section*{Acknowledgements}

The present paper was written in the course of the author's PhD studies under the supervision of Professor Imre Leader. The author is thankful to Prof.\ Leader for his unwavering support.

\end{document}

%% file: bpde-01-intro.tex
\section{Introduction} \label{intro}

Let $d$ be a positive integer and let $\mathcal{A}$ be a finite collection of finite subsets of $\mathbb{Z}^d$. We call each member of $\mathcal{A}$ a \emph{pattern}. Let also $S$ be any subset of $\mathbb{Z}^d$, finite or infinite. We call $S$ a \emph{seed}.

Consider the following process: Throughout, each point of $\mathbb{Z}^d$ is either \emph{occupied} or \emph{vacant}. Initially, all points of the seed are occupied and all other points are vacant. At any moment during the process, we say that a given point \emph{completes} a given subset of $\mathbb{Z}^d$ when it is the unique vacant point in that subset. On each step of the process, all vacant points which complete a translation copy of some pattern in $\mathcal{A}$ become occupied simultaneously. In the absence of such points, the process stops.

We say that the seed $S$ \emph{percolates} if it eventually causes all points of $\mathbb{Z}^d$ to become occupied. We are interested in the following problem: Given a collection of patterns $\mathcal{A}$, to determine the smallest possible size of a percolating seed $S$.

The process we just described is a particular kind of bootstrap percolation. This notion was originally introduced in statistical physics to model phase transitions in magnetic materials. Then, soon thereafter, it was picked up and generalised by cellular automata theorists and combinatorialists. Our specific setup fits best within the broader framework \cite{BSU} of ``$\mathcal{U}$-bootstrap percolation''.

We go on to explore some of the simplest special cases. We work in ascending order of the dimension $d$. For each fixed $d$, we also work in ascending order of the size of $\mathcal{A}$, until we reach a number of patterns where ``most'' collections admit finite percolating seeds. We are going to assume, from now on, that each pattern in $\mathcal{A}$ contains at least two points. (The empty pattern is irrelevant, while putting a single-point pattern in $\mathcal{A}$ trivialises the problem by allowing the empty seed to percolate.)

Consider first the case of $d = 1$ and a single pattern $A$. Then a finite percolating seed certainly does exist. For example, any interval formed out of $\max A - \min A$ consecutive integers does the job. We will show, in Section \ref{prelim}, that this construction actually achieves the smallest size of a percolating seed in the current setting. We view this size as follows: Given any finite subset $A$ of $\mathbb{R}^d$, we write $\lhull A \rhull$ for the convex hull of $A$. Furthermore, given any polytope $F$ in $\mathbb{R}^d$, we write $\vol F$ for the $d$-dimensional volume of $F$.

\begin{theorem} \label{1d} For one pattern $A$ in one dimension, the smallest size of a percolating seed is given by $\vol \lhull A \rhull$. \end{theorem}

We move on to $d = 2$ dimensions next. It is straightforward to see that no single pattern $A$ can admit a finite percolating seed anymore. Indeed, take the outward unit normal $\mathbf{u}$ of any side of $\lhull A \rhull$. For any finite seed $S$, we can find some integer constant $c$ such that $S$ is contained entirely within the half-space $U = \{\mathbf{a} \mid \mathbf{u} \cdot \mathbf{a} \le c\} \subseteq \mathbb{Z}^2$. But then, no matter how large $S$ is, no points of $\mathbb{Z}^2$ outside of $U$ are ever going to become occupied.

What about two patterns $A$ and $B$? If any $\mathbf{u}$ occurs as an outward unit normal in both of $\lhull A \rhull$ and $\lhull B \rhull$, then a finite percolating seed is ruled out by the same argument as before. Following the terminology of \cite{BSU}, we call any such $\mathbf{u}$ a \emph{stable direction} for $A$ and $B$.

Suppose, otherwise, that $A$ and $B$ do not share a stable direction. We are going to need the notion of \emph{mixed volume} -- a certain numerical measure associated with collections of $d$ convex polytopes in $\mathbb{R}^d$. This notion will be reviewed in Section \ref{prelim}.

We are ready now to state our first main result:

\begin{theorem} \label{2d} For two patterns $A$ and $B$ in two dimensions, the smallest size of a percolating seed is given by $\mix(\lhull A \rhull, \lhull B \rhull)$ in the absence of a stable direction. Otherwise, no finite percolating seed exists. \end{theorem}

We give a self-contained, purely combinatorial proof of Theorem \ref{2d} in Section \ref{plane}.

It will be helpful to recast the bootstrap percolation process as a single-player game. The game is played on $\mathbb{Z}^d$, with the initial position given by $S$ and the rules given by~$\mathcal{A}$. On each move, we choose a single vacant point which completes a translation copy of some pattern in $\mathcal{A}$, and we make this point occupied. The goal of the game is to occupy all points of~$\mathbb{Z}^d$. We call this the \emph{bootstrap percolation game} specified by the collection of patterns~$\mathcal{A}$.

Our motivation for the study of this particular flavour of bootstrap percolation comes from the theory of partial difference equations. We proceed now to outline the connections.

Partial difference equations arise when we discretise partial differential equations. We are going to consider a somewhat restricted framework, where a lot of generality has been given up. A more comprehensive treatment may be found, e.g., in the classical textbook \cite{Ch}.

To us, a partial difference equation is specified by a pattern $A$ and a polynomial $P$. The number of points in $A$ must match the number of formal indeterminates in $P$. Let $\eta$ be a $d$-dimensional numerical array; i.e., a function which assigns a number to each point of $\mathbb{Z}^d$. We say that $\eta$ is a \emph{solution} to the partial difference equation specified by $A$ and $P$ when, for each translation copy of $A$, the numbers which $\eta$ assigns to its points cause $P$ to vanish.

Suppose we are given a partial difference equation, and we wish to find a solution to it. One approach would be to look for some clever formula with $d$ integer arguments which we can evaluate at all points of $\mathbb{Z}^d$ simultaneously. This approach may be described as ``global''.

However, a ``local'' approach is also possible, as follows: We examine the points of $\mathbb{Z}^d$ one by one, assigning numbers to them as we go. For each new point under examination, first we check whether it is constrained by the values we have filled in already. If not, we assign a new number to it arbitrarily. Otherwise, if it completes a translation copy of $A$, we must be more careful. Then we make sure that the new assignment causes $P$ to vanish. We are going to assume that in the latter case there always exists a unique valid assignment. This is clearly a strong assumption -- but it is not an unreasonable one, being satisfied in many natural situations.

We can now think of the ``unconstrained'' points as the elements of the seed. They can all be assigned independently and arbitrarily at the beginning. Following this, the ``constrained'' points get assigned one by one according to the rules of the bootstrap percolation game with $\mathcal{A} = \{A\}$. Being able to win the game tells us that the numbers assigned to the seed determine the rest of the solution uniquely. Notice that, in this picture, the polynomial $P$ has been abstracted away completely; the analysis depends solely on the combinatorial properties of the pattern $A$.

One subtlety in the preceding sketch is worth expounding on. Suppose that, on some step during the assignment process, the new point under examination happens to complete two distinct translation copies of $A$ simultaneously. We would get two distinct instances of $P$ which involve that point, and hence also two distinct ways to determine the new number which must be assigned to it. What if they disagree? Then we would run into a contradiction, failing to obtain a solution.

Fortunately, in the setting of a single partial difference equation, this issue is easy to sidestep. We can always find a percolating seed and a winning series of moves from it which avoid such conflicts altogether. The construction will be spelled out in Section \ref{cons-ii}.

However, consistency becomes a much more substantial concern when we turn from individual partial difference equations to systems of them. Then $\mathcal{A}$ comprises multiple patterns, and we must watch out for contradictions whenever the new point under examination completes two or more translation copies of these patterns simultaneously. The issue ceases to be purely combinatorial in this setting -- instead, we obtain consistent systems for some choices of the underlying polynomials, but not for others.

The literature on partial difference equations has already engaged with some questions of consistency; see, for example, \cite{Bo}. On the other hand, the patterns involved are often quite simple from a combinatorial perspective -- such as the four-point pattern given by the vertices of a unit square, or the six-point pattern given by the vertices of a unit-radius regular octahedron.

We develop a formal notion of consistency in Sections \ref{cons-i} and \ref{cons-ii}. Briefly speaking, we say that a system of partial difference equations is consistent if we never run into contradictions as we build solutions to it locally, filling in new values one by one according to the rules of the associated bootstrap percolation game. Our analysis applies to arbitrary collections of patterns~$\mathcal{A}$. Furthermore, when $\mathcal{A}$ admits a finite percolating seed, the smallest such seeds turn out to play a crucial role in the study of consistency.

Our proof of Theorem \ref{2d} yields one consistency criterion for systems of two partial difference equations in two dimensions. The gist is that consistency over a certain small region is necessary and sufficient for consistency on the entire plane. We are going to need the notion of \emph{Minkowski addition} -- a way of summing together convex polytopes in $\mathbb{R}^d$. This notion will be reviewed in Section \ref{prelim}. The terms ``proper'' and ``multilinear'' which appear below signify that the system satisfies certain technical requirements; both of them will be discussed in Section \ref{cons-i}.

We are prepared now to state our second main result:

\begin{theorem} \label{cp} Let $\Gamma$ be a proper multilinear system of two partial difference equations in two dimensions whose patterns $A$ and $B$ admit a finite percolating seed. Then $\Gamma$ is consistent on the plane if and only if it is consistent over the region $\lhull A \rhull + \lhull B \rhull$. \end{theorem} % consistency in the plane

The rest of the paper is structured as follows: Section \ref{prelim} covers the preliminaries. The proof of Theorem \ref{2d} occupies Section \ref{plane}. Sections \ref{cons-i} and \ref{cons-ii} develop the aforementioned formal notion of consistency, concluding with the proof of Theorem \ref{cp}. Section \ref{sd} demonstrates some applications in concrete natural settings. Finally, Section \ref{further} touches upon higher dimensions and outlines a number of potential directions for further research.

%% file: bpde-02-prelim.tex
\section{Preliminaries} \label{prelim}

Fix a collection of patterns $\mathcal{A}$ in $\mathbb{Z}^d$. Of course, replacing any pattern in $\mathcal{A}$ with a translation copy of itself does not in any way affect the bootstrap percolation game specified by $\mathcal{A}$. From now on, we will assume by convention that the members of $\mathcal{A}$ are pairwise distinct modulo translation.

Let $T$ be a subset of $\mathbb{Z}^d$ and let $\mathbf{a}$ be any point in $T$. We define the \emph{degree} of $\mathbf{a}$ in $T$ to be the number of translation copies of some pattern in $\mathcal{A}$ which contain $\mathbf{a}$ and are contained within $T$.

Let $T'$ and $T''$ be two subsets of $\mathbb{Z}^d$ with $T' \subseteq T''$. We define a \emph{combinatorial run} from $T'$ to $T''$ to be any sequence of points in $\mathbb{Z}^d$ which contains each point of $T'' \setminus T'$ exactly once, and which does not contain any other terms besides these points. For simplicity, throughout Sections \ref{prelim} and \ref{plane} we will write just ``run'' instead of ``combinatorial run''. (Section \ref{cons-i} is going to introduce the similar but distinct notion of a ``computational run''.)

A run is \emph{legal} when it constitutes a legal series of moves leading from $T'$ to $T''$ according to the rules of the bootstrap percolation game specified by $\mathcal{A}$. Formally, the run $\mathbf{a}_1$, $\mathbf{a}_2$, $\ldots$, $\mathbf{a}_k$ from $T'$ to $T''$ is legal when the degree of $\mathbf{a}_i$ in $T' \cup \{\mathbf{a}_1, \mathbf{a}_2, \ldots, \mathbf{a}_i\}$ is nonzero for all $i$. We will sometimes consider illegal runs, too. Given an illegal run, we can convert it into a legal one by removing all degree-zero points from it and adjoining them to $T'$ instead.

We go on now to review some basic notions from discrete geometry.

Let $\mathcal{F} = \{F_1, F_2, \ldots, F_k\}$ be a finite collection of convex polytopes in $\mathbb{R}^d$. The \emph{Minkowski sum} of $\mathcal{F}$ is the convex polytope $\{\mathbf{a}_1 + \mathbf{a}_2 + \cdots + \mathbf{a}_k \mid \mathbf{a}_i \in F_i \text{ for all } i\}$. We denote it by $F_1 + F_2 + \cdots + F_k$ or $\operatorname{sum} \mathcal{F}$. Notice that the Minkowski sum of the empty collection $\mathcal{F} = \varnothing$ is the single-point polytope $\mathbf{0}$ given by the origin of our coordinate system for $\mathbb{R}^d$.

Let $f$ be any function which assigns a number to each convex polytope in $\mathbb{R}^d$. We set $\Delta(f \mathbin{|} \mathcal{F}) = \sum_{\mathcal{G} \subseteq \mathcal{F}} (-1)^{|\mathcal{F} \setminus \mathcal{G}|} f(\operatorname{sum} \mathcal{G})$. Then we define the mixed volume of $d$ convex polytopes in $\mathbb{R}^d$ by $\mix(F_1, F_2, \ldots, F_d) = \Delta(\vol \mathbin{|} F_1, F_2, \ldots, F_d)$. So, in particular, the mixed volume and the ordinary volume are related by $\mix(F, F, \ldots, F) = d! \cdot \vol F$.

Suppose that $F$ is an integer convex polytope, meaning that all of its vertices are integer points. We write $\disc F$ for the set of all integer points in $F$. (We consider all of our polytopes to be closed, and so the integer points on the boundary of $F$ are included as well.) The \emph{discrete volume} of $F$ is the number of such points; i.e., $|\disc F|$. We denote it by $\vd F$. It is well-known \cite{Bi} that, for integer convex polytopes, the mixed volume can also be defined in terms of the discrete volume instead of the continuous volume. Or, in other words, $\mix(F_1, F_2, \ldots, F_d) = \Delta(\vd \mathbin{|} F_1, F_2, \ldots, F_d)$. This definition is in fact the more natural one in the context of Theorem \ref{2d}, as we are going to see in Section \ref{plane}.

The mixed volume is well-known \cite{S} to be multi-additive -- meaning that it is additive in each argument with respect to Minkowski addition. The necessary and sufficient condition for $\Delta(f \mathbin{|} F_1, F_2, \ldots, F_k)$ to be multi-additive, considered as a function of $F_1$, $F_2$, $\ldots$, $F_k$, is easily seen to be that $\Delta(f \mathbin{|} F_1, F_2, \ldots, F_{k + 1})$ vanishes identically. We conclude that $\Delta(\vol \mathbin{|} F_1, F_2, \ldots, F_{d + 1}) = 0$ for any $d + 1$ convex polytopes in $\mathbb{R}^d$ and, similarly, $\Delta(\vd \mathbin{|} F_1, F_2, \ldots, F_{d + 1}) = 0$ for any $d + 1$ integer convex polytopes in $\mathbb{R}^d$.

We next provide short proofs for these properties of the mixed volume in two dimensions, so as to put Section \ref{plane} on elementary footing.

Consider any integer convex polygon $F$. Orient the boundary of $F$ counterclockwise, and then partition it into oriented segments by puncturing it at all integer points on it. We denote the multiset of the integer vectors given by these oriented segments by $\side F$. Clearly, for each multiset $V$ of integer vectors whose sum vanishes, there exists an integer convex polygon $F$ (possibly degenerate, but unique modulo translation) such that $\side F = V$. Furthermore, it is straightforward to see that $\side(F_1 + F_2) = \side F_1 \cup \side F_2$.

\begin{lemma} \label{cvdv} In the plane, $\Delta(\vol \mathbin{|} F_1, F_2) = \Delta(\vd \mathbin{|} F_1, F_2)$ for all integer convex polygons $F_1$ and $F_2$. \end{lemma} % continuous volume and discrete volume

\begin{myproof} By Pick's formula together with the union property of the $\side$ operator. \end{myproof}

\begin{lemma} \label{lmv} In the plane, $\Delta(\vol \mathbin{|} F_1, F_2, F_3) = \Delta(\vd \mathbin{|} F_1, F_2, F_3) = 0$ for all integer convex polygons $F_1$, $F_2$, $F_3$. \end{lemma} % linearity of the mixed volume

\begin{myproof} Consider any integer convex polygon $F$ and let $\mathbf{u}_1$, $\mathbf{u}_2$, $\ldots$, $\mathbf{u}_k$ be the elements of $\side F$, ordered counterclockwise. Let also $\mathbf{u}_i \times \mathbf{u}_j$ be the oriented area of the triangle spanned by $\mathbf{u}_i$ and $\mathbf{u}_j$, in this order. By taking a suitable triangulation of $F$, we find that $\vol F = \sum_{1 \le i < j \le k} \mathbf{u}_i \times \mathbf{u}_j$. This formula and the union property of the $\side$ operator together give us $\Delta(\vol \mathbin{|} F_1, F_2, F_3) = 0$. From here, the discrete-volume identity follows as in the proof of Lemma \ref{cvdv}. \end{myproof}

We continue with our formalisation of the notion of a partial difference equation.

For each point $\mathbf{a}$ of $\mathbb{Z}^d$, we introduce a new formal indeterminate $x_\mathbf{a}$. Let $X$ be the set of all such indeterminates. Fix any field $\mathbb{F}$, and consider the polynomial ring $\mathbb{F}[X]$.

Let $P$ be any polynomial in $\mathbb{F}[X]$. We define the \emph{pattern} of $P$ to be the set of all points $\mathbf{a}$ in $\mathbb{Z}^d$ such that $x_\mathbf{a}$ occurs in $P$. We consider translation by $\mathbf{t}$, with $\mathbf{t} \in \mathbb{Z}^d$, to map $P$ onto the element of $\mathbb{F}[X]$ obtained from $P$ by simultaneously, for each $\mathbf{a} \in \mathbb{Z}^d$, substituting every occurrence of $x_\mathbf{a}$ in $P$ with an occurrence of $x_{\mathbf{a} + \mathbf{t}}$.

We define the \emph{partial difference equation} specified by $P$ to be the system of polynomial equations over the unknowns $X$ formed by all translation copies of $P$. Fix any field $\mathbb{G}$ which extends $\mathbb{F}$ and let $\eta$ be any function which maps each point $\mathbf{a}$ of $\mathbb{Z}^d$ onto an element $\eta_\mathbf{a}$ of $\mathbb{G}$. We write $P(\eta)$ for the element of $\mathbb{G}$ obtained when $P$ is evaluated with $x_\mathbf{a} = \eta_\mathbf{a}$ for all points $\mathbf{a}$~in~$\mathbb{Z}^d$. We define $\eta$ to be a \emph{solution} to the partial difference equation specified by $P$ when $Q(\eta) = 0$ for all translation copies $Q$ of $P$.

We may also consider partial difference equations over smaller regions within $\mathbb{Z}^d$. Let $T$ be any subset of $\mathbb{Z}^d$. We define the partial difference equation specified by $P$ over $T$ to be the system of polynomial equations over the unknowns $X_T = \{x_\mathbf{a} \mid \mathbf{a} \in T\}$ formed by all translation copies of $P$ whose patterns are contained within $T$. Consider now any function $\eta$ which maps each point $\mathbf{a}$ of $T$ onto an element $\eta_\mathbf{a}$ of $\mathbb{G}$. We define $\eta$ to be a solution to the partial difference equation specified by $P$ over $T$ when $Q(\eta) = 0$ for all translation copies $Q$ of $P$ such that the pattern of $Q$ is contained within $T$.

Let $\mathcal{P}$ be a finite collection of polynomials in $\mathbb{F}[X]$. We define the \emph{system} of partial difference equations $\Gamma$ specified by $\mathcal{P}$ to be the union of the systems of polynomial equations over the unknowns $X$ given by the partial difference equations specified by the members of $\mathcal{P}$. We also define a solution to $\Gamma$ to be any solution shared by these partial difference equations. The restrictions of these notions to smaller regions $T$ within $\mathbb{Z}^d$ are handled as before.

We write $\sol \Gamma$ for the set of all solutions to $\Gamma$, and we also write $\sol_T \Gamma$ for the set of all solutions to $\Gamma$ over the region $T$. Notice that, in both of these notations, the field $\mathbb{G}$ out of which we are drawing the values in the solutions is understood implicitly.

We conclude with some preliminary observations tying partial difference equations to the bootstrap percolation game.

Consider any system $\Gamma$ of partial difference equations whose specifying collection $\mathcal{P}$ consists entirely of homogeneous linear polynomials. Of course, in this setting, the solutions to $\Gamma$ form a linear space over $\mathbb{G}$. By a slight abuse of notation, we are going to write $\sol \Gamma$ for this linear space as well.

Given any pattern $A$, let $\add A$ be the partial difference equation specified by the polynomial $\sum_{\mathbf{a} \in A} x_\mathbf{a}$. Similarly, given any finite collection of patterns $\mathcal{A}$, let $\add \mathcal{A}$ be the system comprised of the partial difference equations $\add A$ over all members $A$ of $\mathcal{A}$.

\begin{lemma} \label{dbps} Suppose that $S$ is a percolating seed for $\mathcal{A}$. Then $|S| \ge \dim \sol \add \mathcal{A}$. \end{lemma} % dimension bound for percolating seeds

\begin{myproof} Because any solution $\eta$ to $\add \mathcal{A}$ is uniquely determined by its restriction to $S$. Indeed, if we know what numbers $\eta$ assigns to the points of $S$, then by working our way outwards from $S$, as per any winning series of moves in the bootstrap percolation game specified by $\mathcal{A}$, one by one we can figure out the numbers which $\eta$ assigns to all other points of $\mathbb{Z}^d$ as well. \end{myproof}

We can now establish Theorem \ref{1d} with no difficulty whatsoever.

\begin{myproof}[Proof of Theorem \ref{1d}] Let $s = \vol \lhull A \rhull$ and $S_{\Join} = \{1, 2, \ldots, s\}$. Consider the run $R$ given by $0$, $s + 1$, $-1$, $s + 2$, $-2$, $s + 3$, $\ldots$, from $S_{\Join}$ to $\mathbb{Z}$. Each point in it is of unit degree exactly, and so the seed $S_{\Join}$ percolates. This takes care of the upper bound.

We continue with the lower bound. By Lemma \ref{dbps}, it suffices to show that $\dim \sol \add A \ge |S_{\Join}|$. This, in turn, will follow if we can show that every assignment of values to the points of $S_{\Join}$ can be extended so as to become a solution to $\add A$ over $\mathbb{Z}$. However, the latter objective is straightforward: We can fill in the new values one by one, working our way outwards from $S_{\Join}$, as per the winning series of moves $R$. \end{myproof}

Our proof for the lower bound of Theorem \ref{2d} in Section \ref{plane} is going to follow a similar high-level strategy, albeit complicated by issues of consistency.

%% file: bpde-03-plane.tex
\section{The Plane} \label{plane}

Here, we prove Theorem \ref{2d}. Fix, throughout this section, $\mathcal{A} = \{A, B\}$ with $A$ and $B$ being two finite subsets of $\mathbb{Z}^2$, each containing two or more points, such that $A$ and $B$ do not share a stable direction.

We begin with some purely discrete-geometric observations:

\begin{lemma} \label{vs} Let $F$ be an integer convex polygon in $\mathbb{Z}^2$. Then each vertex of $F$ is of degree at most $2$ in $\disc F$. Furthermore, if all vertices of $F$ are of degree exactly $2$ in $\disc F$, then $F$ is of the form $\lhull A \rhull + \lhull B \rhull + G$ for some integer convex polygon $G$. \end{lemma} % vertices

\begin{myproof} The first part is clear as each vertex of $F$ can belong to at most one translation copy of $\lhull A \rhull$ contained within $F$, and similarly for $\lhull B \rhull$. Suppose, next, that each vertex of $F$ does belong to a translation copy of $\lhull A \rhull$ and a translation copy of $\lhull B \rhull$ both of which are contained within $F$.

Choose any side $u$ of $\lhull A \rhull$ and consider the support line $\ell$ of $F$ parallel to $u$ such that $F$ lies on the same side of $\ell$ as $\lhull A \rhull$ does relative to the line through $u$. Then $\ell$ cannot touch $F$ at a single vertex, as it would have been impossible for that vertex to belong to a translation copy of $\lhull A \rhull$ contained within $F$. So $\ell$ must touch $F$ along a side $v$. By a similar argument, we find that also $|u| \le |v|$.

The same reasoning applies to each side of $\lhull B \rhull$ as well. Since $A$ and $B$ do not share a stable direction, we get that $\side \lhull A \rhull$ and $\side \lhull B \rhull$ are disjoint and their multiset union is contained within the multiset $\side F$. This gives us $G$ by means of $\side G = \side F \setminus (\side \lhull A \rhull \cup \side \lhull B \rhull)$. \end{myproof}

For convenience, we might speak of a run from, or to, some integer convex polygon $F$ when we really mean a run from, or to, the set $\disc F$.

Both Lemmas \ref{zs} and \ref{se} below are about the existence of runs with certain desirable properties. We construct these runs ``in reverse'' -- by taking points away from the final set instead of adding points to the initial set. On each step, the point being taken away is a vertex of the convex hull of the remaining points, chosen out of all such vertices according to considerations which will be specified shortly.

We say that a point in a given run is \emph{safe} when it completes a translation copy of $A + B$ in which it serves as a vertex of the convex hull.

\begin{lemma} \label{zs} There exists a run from $\varnothing$ to $\lhull A \rhull + \lhull B \rhull$ where: (i) The number of points of degree $0$ is $\mix(\lhull A \rhull, \lhull B \rhull)$; (ii) The unique point of degree $2$ is safe; and (iii) There are no points of degree $3$ or more. \end{lemma} % zero to sum

\begin{myproof} We obtain the desired run by boiling $\lhull A \rhull + \lhull B \rhull$ down to $\varnothing$, as described above. By Lemma \ref{vs}, on each step except for the first one we can take away a point of degree either $0$ or $1$ in the set of all remaining points. This ensures (ii) and (iii).

For each $i$, let $\lambda_i$ be the number of points of degree $i$ in the run just constructed. Then $\lambda_2 = 1$ and $\lambda_0 + \lambda_1 + \lambda_2 = |\disc(\lhull A \rhull + \lhull B \rhull)|$ by the preceding discussion. Furthermore, the number of integer translation copies of $A$ contained within $\lhull A \rhull + \lhull B \rhull$ is $|\disc \lhull B \rhull|$; the analogous number for $B$ is $|\disc \lhull A \rhull|$; and $0 \cdot \lambda_0 + 1 \cdot \lambda_1 + 2 \cdot \lambda_2$ amounts to the sum of these two quantities.

From here, we find that $\lambda_0 = \vd(\lhull A \rhull + \lhull B \rhull) - \vd \lhull A \rhull - \vd \lhull B \rhull + 1$. This settles (i), too, as the right-hand side equals $\mix(\lhull A \rhull, \lhull B \rhull)$ by Lemma \ref{cvdv}. \end{myproof}

\begin{lemma} \label{se} There exists a run from $\lhull A \rhull + \lhull B \rhull$ to $\mathbb{Z}^2$ where: (i) There are no points of degree~$0$; (ii) Each point of degree $2$ is safe; and (iii) There are no points of degree $3$ or more. \end{lemma} % sum to everything

\begin{myproof} Fix any run $\mathbf{g}_1 = \mathbf{0}$, $\mathbf{g}_2$, $\mathbf{g}_3$, $\ldots$ from $\varnothing$ to $\mathbb{Z}^2$ such that, for all $i$, the set $\{\mathbf{g}_1, \mathbf{g}_2, \ldots, \mathbf{g}_i\}$ coincides with the image under the $\disc$ operator of its own convex hull $G_i$. Let also $F_i = \lhull A \rhull + \lhull B \rhull + G_i$.

Consider any $V$ with $\disc F_i \subsetneq V \subsetneq \disc F_{i + 1}$. By Lemma \ref{vs}, there exists a vertex $\mathbf{v}$ of $\lhull V \rhull$ with degree either $0$ or $1$ in $V$. But $\mathbf{v}$ cannot be a vertex of $F_i$, as all vertices of $F_i$ are already of degree $2$ in $\disc F_i$. So we can boil $F_{i + 1}$ down to $F_i$, as described above, in order to obtain a run $R_i$ from $F_i$ to $F_{i + 1}$ where all points except for the final one are of degree either $0$ or $1$; while the final point is of degree $2$ and safe.

Fix $i$ and consider the run $R$ from $\varnothing$ to $F_i$ obtained by concatenating the run of Lemma~\ref{zs} together with the runs $R_1$, $R_2$, $\ldots$, $R_{i - 1}$ just constructed. For each $j$, let $\lambda_j$ be the number of points of degree $j$ in $R$. Then $\lambda_0 + \lambda_1 + \lambda_2 = |\disc F_i|$ and $\lambda_2 = |\disc G_i|$ by the preceding discussion. Furthermore, the number of integer translation copies of $A$ contained within $F_i$ is $|\disc(\lhull B \rhull + G_i)|$; the analogous number for $B$ is $|\disc(\lhull A \rhull + G_i)|$; and $0 \cdot \lambda_0 + 1 \cdot \lambda_1 + 2 \cdot \lambda_2$ works out to the sum of these two quantities.

From here, we find that $\lambda_0 = \vd(\lhull A \rhull + \lhull B \rhull + G_i) - \vd(\lhull A \rhull + G_i) - \vd(\lhull B \rhull + G_i) + \vd G_i$. The right-hand side equals $\mix(\lhull A \rhull, \lhull B \rhull)$ by Lemmas \ref{cvdv} and \ref{lmv}. We conclude, in light of Lemma \ref{zs}, that none of the runs $R_1$, $R_2$, $R_3$, $\ldots$ contain any points of degree zero. Thus the concatenation of these runs satisfies all three of (i), (ii), and (iii). \end{myproof}

We next address the problem of consistency:

\begin{lemma} \label{lc} Let $\mathbf{v}$ be a vertex of $\lhull A \rhull + \lhull B \rhull$. Suppose that $\eta$ is a solution to $\add \mathcal{A}$ over $(A + B) \setminus \{\mathbf{v}\}$. Then we can extend $\eta$, by assigning a value to $\mathbf{v}$ as well, so as to obtain a solution to $\add \mathcal{A}$ over $A + B$. \end{lemma} % linear consistency

\begin{myproof} Let $P$ be the specifying polynomial of $\add A$. For each $\mathbf{b} \in B$, take the image of $P$ upon translation by $\mathbf{b}$, and let these translation copies form the collection of polynomials $\mathcal{P}$. The pattern of each member of $\mathcal{P}$ is contained within $A + B$; and one of the members of $\mathcal{P}$ is the unique translation copy of $P$ such that its pattern contains $\mathbf{v}$ and is contained within~$A + B$.

Define the polynomial $Q$ and the collection of polynomials $\mathcal{Q}$ similarly, with $A$ and $B$ swapped. The desired result now follows because the sum of all members of $\mathcal{P}$ equals the sum of all members of $\mathcal{Q}$. Indeed, both of these sums work out to $\sum x_{\mathbf{a} + \mathbf{b}}$ over all $\mathbf{a} \in A$ and $\mathbf{b} \in B$. \end{myproof}

What remains is merely to put the pieces together:

\begin{myproof}[Proof of Theorem \ref{2d}] By Lemmas \ref{zs} and \ref{se}, we get a run from $\varnothing$ to $\mathbb{Z}^2$ with $\mix(\lhull A \rhull, \lhull B \rhull)$ points of degree zero. The set $S_{\Join}$ formed by these points is a percolating seed for $\mathcal{A}$. This settles the upper bound.

We continue with the lower bound. By Lemma \ref{dbps}, it suffices to show that $\dim \sol \add \mathcal{A} \ge |S_{\Join}|$. This, in turn, will follow if we can show that every assignment of values to the points of $S_{\Join}$ can be extended so as to become a solution to $\add \mathcal{A}$ over $\mathbb{Z}^2$. We fill in the new values one by one, working our way outwards from $S_{\Join}$, as per the winning series of moves constructed in Lemmas \ref{zs} and \ref{se}. It is only the moves of degree $2$ which might conceivably cause a contradiction. However, since all such moves are safe, Lemma \ref{lc} guarantees that in fact no contradictions are going to arise. \end{myproof}

%% file: bpde-04-cons-i.tex
\section{Consistency I} \label{cons-i}

Let $\Gamma$ be a system of partial difference equations specified by the collection of polynomials $\mathcal{P} \subseteq \mathbb{F}[X]$. We study the solutions of $\Gamma$ over the region $T \subseteq \mathbb{Z}^d$, with values drawn out of the field $\mathbb{G}$. We assume throughout Sections \ref{cons-i} and \ref{cons-ii} that the patterns of the polynomials which specify $\Gamma$ are pairwise distinct upon translation, and we denote the collection of these patterns by $\mathcal{A}$.

We aim to solve $\Gamma$ ``locally'', filling in values one by one. Consider any translation copy $P$ of some member of $\mathcal{P}$ whose pattern $A$ is contained within $T$, and let $\mathbf{a}$ be any point of $A$. Suppose that we have already assigned a value to each point of $A \setminus \{\mathbf{a}\}$, by means of the assignment~$\eta_\text{Part}$. We wish, as announced in the introduction, to focus on the setting where the value to be assigned to $\mathbf{a}$ can be deduced uniquely from the equation $P = 0$.

Thus we stipulate that $\Gamma$ must be multilinear; i.e., that each member of $\mathcal{P}$ must be linear in every one of its indeterminates. It is only multilinear systems that we are going to consider in Sections \ref{cons-i} and \ref{cons-ii}. (One generalisation of our framework which allows the multilinearity constraint to be relaxed somewhat will be touched upon briefly in Section \ref{further}.)

We can now write $P = P_\text{Num} + x_\mathbf{a}P_\text{Denom}$, with both of $P_\text{Num}$ and $P_\text{Denom}$ in $\mathbb{F}[X \setminus \{x_\mathbf{a}\}]$. We call these polynomials the \emph{numerator} and \emph{denominator} of $P$ at $\mathbf{a}$. The new value to be assigned to $\mathbf{a}$ is then given by $-P_\text{Num}(\eta_\text{Part})/P_\text{Denom}(\eta_\text{Part})$, provided that $P_\text{Denom}(\eta_\text{Part}) \neq 0$.

We next formalise the notion of filling in new values one by one. There is one minor technical concession we need to make: We are going to allow these values to be recomputed multiple times, in different ways, so as to account for the possibility that $\Gamma$ might be inconsistent. We are also going to focus exclusively on the setting where a finite percolating seed exists. Most of the material in Sections \ref{cons-i} and \ref{cons-ii} is straightforward to adapt with infinite percolating seeds as well -- though Theorems \ref{cp} and \ref{css} constitute two significant exceptions.

Let $S$ be any finite percolating seed for $T$ in the bootstrap percolation game specified by $\mathcal{A}$. We define a \emph{computational run} from $S$ to $T$ to be any ordered pair $R = (\alpha, \Psi)$ of two sequences $\alpha$ and $\Psi$ such that: The sequence $\alpha$ is indexed from $1$; its first $s = |S|$ elements are the points of $S$; each point of $T$ occurs in it, with repetitions allowed; and all of its elements are points of $T$. The sequence $\Psi$ is indexed from $s + 1$; for each index $i$, its $i$-th term $\Psi(i)$ is a subset of $\{1, 2, \ldots, i\}$ which contains $i$; as the index $j$ ranges over $\Psi(i)$, the points $\alpha(j)$ are pairwise distinct; and these points form a translation copy of some pattern in $\mathcal{A}$.

The intuition is that $\alpha$ tells us what order to visit the points of $T$ in, while $\Psi$ tells us which equations in particular to employ when we are computing the new values to be assigned to these points. We make this intuition precise below.

Let $\eta_\text{Seed}$ be an assignment which maps each point of $S$ onto a value in $\mathbb{G}$. We define the \emph{computation} from $\eta_\text{Seed}$ based on the computational run $R$ to be the sequence $\xi$ in $\mathbb{G}$ with the same indexing interval as $\alpha$ whose $i$-th term $\xi(i)$ is constructed as follows: For each $1 \le i \le s$, simply $\xi(i) = \eta_\text{Seed}(\alpha(i))$. For each $i \ge s + 1$, let $P$ be the unique translation copy of some member of $\mathcal{P}$ whose pattern is given by $\{\alpha(j) \mid j \in \Psi(i)\}$, with numerator $P_\text{Num}$ and denominator $P_\text{Denom}$ at $\alpha(i)$. Let also $\eta_\text{Step}$ be the assignment given by $\eta_\text{Step}(\alpha(j)) = \xi(j)$ as the index $j$ ranges over $\Psi(i) \setminus \{i\}$. Then $\xi(i) = -P_\text{Num}(\eta_\text{Step})/P_\text{Denom}(\eta_\text{Step})$, provided that $P_\text{Denom}(\eta_\text{Step}) \neq 0$. Otherwise, if $P_\text{Denom}(\eta_\text{Step}) = 0$ for some index $i \ge s + 1$, the computation from $\eta_\text{Seed}$ based on the computational run $R$ is not well-defined.

For simplicity, throughout Sections \ref{cons-i} and \ref{cons-ii} we will write just ``run'' instead of ``computational run''. (Recall that Section \ref{prelim} introduced the similar but distinct notion of a ``combinatorial run''.)

Clearly, if the computation $\xi$ from $\eta_\text{Seed}$ based on $R$ is well-defined, then $\Gamma$ admits at most one solution over $T$ whose restriction to $S$ coincides with $\eta_\text{Seed}$. Furthermore, if $\Gamma$ does in fact admit such a solution $\eta$, then for all $i$ in the indexing interval of $\alpha$ and $\xi$ it must hold that $\xi(i) = \eta(\alpha(i))$.

We proceed now to work our way around the issue of potential division by zero in our computations. Let $x_1$, $x_2$, $\ldots$, $x_s$ be new formal indeterminates, and set $\mathbb{G} = \mathbb{F}(x_1, x_2, \ldots, x_s)$. Let also $\eta_\text{Var}$ be any bijection between $S$ and $\{x_1, x_2, \ldots, x_s\}$. We define the \emph{free computation} associated with $R$ to be the computation from $\eta_\text{Var}$ based on $R$. The intuition is that the free computation ``captures'' all other computations based on $R$, as we should be able to obtain these other computations from it by substituting suitable values into the indeterminates.

We define $\Gamma$ to be \emph{proper} over $T$ when, for each finite percolating seed $S$ and every run $R$ from $S$ to $T$, the free computation based on $R$ is well-defined. We restrict our considerations from now on to proper systems because, roughly speaking, for them division by zero is not a major problem.

It might seem at first glance that propriety will usually be difficult to verify for concrete $\Gamma$. However, in practice, we can often verify it without too much trouble in the settings we care about. For example, if $\Gamma$ is linear, then all denominators will be nonzero constants. Otherwise, in nonlinear settings, one potential strategy for the verification of propriety would be to look for a suitable ``propriety certificate'', as outlined below.

We define the solution $\eta$ to $\Gamma$ over $T$ to be \emph{proper} when, for each translation copy $P$ of some member of $\mathcal{P}$ whose pattern $A$ is contained within $T$, and for each point $\mathbf{a}$ of $A$, the denominator $P_\text{Denom}$ of $P$ at $\mathbf{a}$ satisfies $P_\text{Denom}(\eta) \neq 0$. Clearly, if $\eta$ is a proper solution to $\Gamma$ over $T$ and $\eta_\text{Seed}$ is the restriction of $\eta$ to an arbitrary finite percolating seed $S$, then the computation $\xi$ from $\eta_\text{Seed}$ based on $R$ will be well-defined for every run $R$ from $S$ to $T$.

On the other hand, a rational function which attains any nonzero value whatsoever cannot possibly be the zero rational function. Thus we can certify propriety by exhibiting any proper solution to $\Gamma$ over $T$. Furthermore, if all members of $\mathcal{P}$ are integer-coefficient polynomials, then it suffices to exhibit a proper solution to $\Gamma$ over $T$ modulo any fixed prime $p$. Or, in other words, a proper solution to $\Gamma$ over $T$ where the values are drawn out of the finite field $\mathbb{F}_p$.

%% file: bpde-05-cons-ii.tex
\section{Consistency II} \label{cons-ii}

We assume throughout this section that $\Gamma$ is multilinear and proper over $T$. We define $\Gamma$ to be \emph{consistent} from $S$ to $T$ when, for every run $R = (\alpha, \Psi)$ from $S$ to $T$, the free computation $\xi$ associated with $R$ assigns identical values to identical points, so that $\xi(i) = \xi(j)$ whenever $\alpha(i) = \alpha(j)$.

It is straightforward to see that $\Gamma$ is consistent from $S$ to $T$ if and only if $\Gamma$ admits a solution over $T$ with values drawn out of the field $\mathbb{F}(x_1, x_2, \ldots, x_s)$ where the points of $S$ are in a bijective correspondence with the indeterminates $x_1$, $x_2$, $\ldots$, $x_s$.

Suppose that we are given concrete $\Gamma$, $S$, $T$ and we wish to verify the consistency of $\Gamma$ from $S$ to $T$ computationally. The most natural approach would be as follows: First we fix any combinatorial run $R_\text{Comb}$ from $S$ to $T$. Then we go through $R_\text{Comb}$, computing the values to be assigned to the points of $T \setminus S$ one by one. On each move, if the point $\mathbf{a}$ under consideration is of degree $e$ in $R_\text{Comb}$, we compute the value to be assigned to it in all of the $e$ distinct ways afforded by $R_\text{Comb}$. If these candidate values do not all coincide, we declare inconsistency. Otherwise, if they do all coincide, we assign that value to $\mathbf{a}$ and we continue onwards with $R_\text{Comb}$. Once all moves of $R_\text{Comb}$ have been made successfully, we declare $\Gamma$ to be consistent from $S$ to $T$.

We go on now to show that the smallest percolating seeds play a crucial role in the study of consistency.

\begin{theorem} \label{css} Let $\Gamma$ be a multilinear system of partial difference equations proper over the region $T$. Let also $S$ be a finite percolating seed for $T$. Suppose that $\Gamma$ is consistent from $S$ to~$T$. Then $s = |S|$ is the smallest size of a percolating seed for $T$. Furthermore, if $S^\star$ is any other finite percolating seed for $T$ of the same size $s$, then $\Gamma$ is also consistent from $S^\star$ to $T$. \end{theorem} % consistency and smallest seeds

Observe that Theorem \ref{css} allows us to deal away with the seed in the definition of consistency. Or, in other words, it reveals consistency to be a property solely of the system and the region. We define $\Gamma$ to be consistent over $T$ when there exists any finite percolating seed $S$ for $T$ such that $\Gamma$ is consistent from $S$ to $T$. By Theorem \ref{css}, if $\Gamma$ is consistent over $T$, then it is also consistent from $S$ to $T$ whenever $S$ is a smallest percolating seed for $T$.

For convenience, we might also speak of the consistency of $\Gamma$ over some integer convex polygon $F$ when we really mean the consistency of $\Gamma$ over the region $\disc F$.

For the proof of Theorem \ref{css}, we need one notion from abstract algebra. Given a field $\mathbb{F}$ and an extension of it $\mathbb{G}$, the \emph{transcendence degree} of $\mathbb{G}$ over $\mathbb{F}$ is the greatest number of elements of $\mathbb{G}$ which are algebraically independent over $\mathbb{F}$. The basic properties of transcendence degrees are covered, for example, in the classical textbook \cite{DF}. The abstract-algebraic apparatus we require is contained within Lemmas \ref{gf} and \ref{fg} below.

Let $\mathbf{f} = (f_1, f_2, \ldots, f_k)$ be an ordered tuple of rational functions, each with $\ell$ arguments and coefficients drawn out of $\mathbb{F}$. Similarly, let $\mathbf{g} = (g_1, g_2, \ldots, g_\ell)$ be an ordered tuple of rational functions, each with $k$ arguments and coefficients drawn out of $\mathbb{F}$. Given any $\mathbf{u} \in \mathbb{G}^\ell$, we write $\mathbf{f}(\mathbf{u})$ for $(f_1(\mathbf{u}), f_2(\mathbf{u}), \ldots, f_k(\mathbf{u}))$; and analogously for $\mathbf{g}(\mathbf{v})$ with $\mathbf{v} \in \mathbb{G}^k$. Finally, let $\mathbf{x} = (x_1, x_2, \ldots, x_\ell)$ and $\mathbf{y} = (y_1, y_2, \ldots, y_k)$ be two ordered tuples of new formal indeterminates.

\begin{lemma} \label{gf} Suppose that $\mathbf{g}(\mathbf{f}(\mathbf{x})) = \mathbf{x}$. Then $k \ge \ell$. \end{lemma}

\begin{myproof} Consider the field $\mathbb{F}(f_1(\mathbf{x}), f_2(\mathbf{x}), \ldots, f_k(\mathbf{x}))$. On the one hand, it is a subfield of $\mathbb{F}(x_1, x_2, \ldots, x_\ell)$. On the other hand, by virtue of $g_i(f_1(\mathbf{x}), f_2(\mathbf{x}), \ldots, f_k(\mathbf{x})) = x_i$, it must contain the rational functions $x_1$, $x_2$, $\ldots$, $x_\ell$, and so it is also an extension of $\mathbb{F}(x_1, x_2, \ldots, x_\ell)$. We conclude that these two fields coincide. However, the transcendence degree of $\mathbb{F}(f_1(\mathbf{x}), f_2(\mathbf{x}),\allowbreak \ldots, f_k(\mathbf{x}))$ over $\mathbb{F}$ cannot exceed $k$; whereas the transcendence degree of $\mathbb{F}(x_1, x_2, \ldots, x_\ell)$ over $\mathbb{F}$ equals $\ell$ exactly. \end{myproof}

\begin{lemma} \label{fg} Suppose that $\mathbf{g}(\mathbf{f}(\mathbf{x})) = \mathbf{x}$ and $k = \ell$. Then also $\mathbf{f}(\mathbf{g}(\mathbf{y})) = \mathbf{y}$. \end{lemma}

\begin{myproof} As in the proof of Lemma \ref{gf}, we get that the transcendence degree of $\mathbb{F}(f_1(\mathbf{x}), f_2(\mathbf{x}),\allowbreak \ldots, f_k(\mathbf{x}))$ over $\mathbb{F}$ equals $k$ exactly, and so $f_1(\mathbf{x})$, $f_2(\mathbf{x})$, $\ldots$, $f_k(\mathbf{x})$ must be algebraically independent over $\mathbb{F}$. Hence, there exists an isomorphism between the fields $\mathbb{F}(f_1(\mathbf{x}), f_2(\mathbf{x}), \ldots, f_k(\mathbf{x}))$ and $\mathbb{F}(y_1, y_2,\allowbreak \ldots, y_k)$ which maps $f_i(\mathbf{x})$ onto $y_i$ for all $i$.

We next apply $\mathbf{f}$ to both sides of $\mathbf{g}(\mathbf{f}(\mathbf{x})) = \mathbf{x}$ so as to obtain that $\mathbf{f}(\mathbf{g}(\mathbf{f}(\mathbf{x}))) = \mathbf{f}(\mathbf{x})$ holds in $\mathbb{F}(f_1(\mathbf{x}), f_2(\mathbf{x}), \ldots, f_k(\mathbf{x}))$. By means of the aforementioned isomorphism, it follows that also $\mathbf{f}(\mathbf{g}(\mathbf{y})) = \mathbf{y}$ holds in $\mathbb{F}(y_1, y_2, \ldots, y_k)$, as desired. \end{myproof}

We are ready now to establish Theorem \ref{css}.

\begin{myproof}[Proof of Theorem \ref{css}] First let $S^\star$ be any finite percolating seed for $T$, with $t = |S^\star|$. Fix any run $R$ from $S$ to $T$ and let the free computation based on $R$ with indeterminates $x_1$, $x_2$, $\ldots$, $x_s$ assign the rational functions $f_1$, $f_2$, $\ldots$, $f_t$ to the points of $S^\star$. Conversely, fix also any run $R^\star$ from $S^\star$ to $T$, without repetitions, and let the free computation based on $R^\star$ with indeterminates $y_1$, $y_2$, $\ldots$, $y_t$ assign the rational functions $g_1$, $g_2$, $\ldots$, $g_s$ to the points of $S$.

Since $\Gamma$ is consistent from $S$ to $T$, by ``chaining together'' $R$ and $R^\star$ we get that $\mathbf{g}(\mathbf{f}(\mathbf{x})) = \mathbf{x}$. By Lemma \ref{gf}, it follows that $s \le t$. Hence, $s$ is indeed the smallest size of a percolating seed for~$T$.

Suppose, from now on, that $s = t$. Consider the solution $\eta$ of $\Gamma$ over $T$ where the points of $S$ are assigned the indeterminates $x_1$, $x_2$, $\ldots$, $x_s$. Since $\Gamma$ is proper over $T$, by chaining together $R^\star$ and $R$ we find that substituting each $x_i$ with $g_i$ in the values of $\eta$ yields well-defined rational functions. These rational functions form a new solution $\eta^\star$ to $\Gamma$ over $T$. Furthermore, since $\mathbf{f}(\mathbf{g}(\mathbf{y})) = \mathbf{y}$ by Lemma \ref{fg}, in $\eta^\star$ the points of $S^\star$ must be assigned the indeterminates $y_1$, $y_2$, $\ldots$,~$y_t$. Thus $\Gamma$ is indeed consistent from $S^\star$ to $T$ as well. \end{myproof}

Our notion of consistency is ``combinatorially fragile'', and one ought to be careful when working with it. One crucial subtlety is that consistency over a region does not automatically imply consistency over its subregions. The intuition is that a subregion, despite being smaller, might nevertheless be more complicated from the perspective of the bootstrap percolation game.

For example, in $d = 2$ dimensions, consider the system $\Gamma$ specified by the collection of polynomials $\mathcal{P} = \{x_{(0, 0)} - x_{(1, 0)},\; x_{(0, 0)} + 2x_{(2, 0)} - x_{(0, 1)},\; 2x_{(0, 0)} + x_{(2, 0)} - x_{(2, 1)}\}$. It is consistent over the region $T = \{0, 1, 2\} \times \{0, 1\}$ but not over its subregion $T \setminus \{(1, 0)\}$.

Still, if $\Gamma$ is consistent from $S$ to $T$, and $T^\star$ is a subregion of $T$ such that $S^\star = S \cap T^\star$ is a finite percolating seed for $T^\star$, then of course $\Gamma$ will be consistent from $S^\star$ to $T^\star$ as well.

When $\Gamma$ is comprised of a single partial difference equation, it exhibits consistency for more or less trivial reasons. The key combinatorial observation is as follows:

\begin{lemma} \label{es} Suppose that $\mathcal{A} = \{A\}$. Then there exists a combinatorial run from $\varnothing$ to $\mathbb{Z}^d$ where all points are of degree either $0$ or $1$. \end{lemma} % single-equation system

\begin{myproof} Fix a basis $U$ of $\mathbb{R}^d$ such that, for each $\mathbf{u} \in U$, both of the supporting hyperplanes of $A$ perpendicular to $\mathbf{u}$ touch $A$ at a single point. It is straightforward to construct an infinite sequence of nested boxes $F_1 \subseteq F_2 \subseteq F_3 \subseteq \cdots$ in $\mathbb{R}^d$ such that $F_1 \cap \mathbb{Z}^d = \{\mathbf{0}\}$; $\bigcup_{i \in \mathbb{N}} F_i = \mathbb{R}^d$; for all $i$, the vectors $\{\pm \mathbf{u} \mid \mathbf{u} \in U\}$ serve as outward normals for the faces of $F_i$; and, for all $i$, the set $(F_{i + 1} \setminus F_i) \cap \mathbb{Z}^d$ is contained within some face of $F_{i + 1}$. We now let our desired combinatorial run visit the points of $\mathbb{Z}^d$ in the same order as they appear inside of these boxes. \end{myproof}

\begin{theorem} \label{ce} Let $\Phi$ be a multilinear partial difference equation proper over the finite region~$T$. Then $\Phi$ is consistent over $T$. \end{theorem} % consistency of a single equation

It is not too difficult to adapt this result to infinite regions $T$, too; however, if the percolating seed implicit in the proof turns out to be infinite, then we cannot make use of the inherently finitary Theorem \ref{css} anymore, and instead of consistency over $T$ we must speak only of consistency from that percolating seed to $T$.

\begin{myproof} Let $R_\text{Comb}$ be the sub-run of the combinatorial run constructed in Lemma \ref{es} comprised of the points of $T$. We now simply assign values to the points of $T$ one by one as per $R_\text{Comb}$. \end{myproof}

We move on next to two partial difference equations in two dimensions, as addressed by Theorem \ref{cp}.

\begin{myproof}[Proof of Theorem \ref{cp}] The case when $|A| = 1$ or $|B| = 1$ is clear, and so we assume from now on that $|A| \ge 2$ and $|B| \ge 2$ as in Section \ref{plane}. Let $R_\text{Win}$ be the combinatorial run from $\varnothing$ to $\mathbb{Z}^2$ constructed in the proof of Theorem \ref{2d}. Since the seed $S_{\Join}$ given by $R_\text{Win}$ is percolating for both of $\lhull A \rhull + \lhull B \rhull$ and $\mathbb{Z}^2$, necessity follows immediately.

We continue with sufficiency. To the points of $S_{\Join}$, we assign the new formal indeterminates $x_1$, $x_2$, $\ldots$, $x_s$ with $s = \mix(\lhull A \rhull, \lhull B \rhull)$. Then we fill in the remaining values on the plane one by one, working our way through $R_\text{Win}$. Consider any move $\mathbf{v}$ in $R_\text{Win}$ of degree $2$. Suppose that there has been no contradiction so far. Our goal is to show that no contradiction will arise at $\mathbf{v}$, either.

Since $\mathbf{v}$ is safe, it completes a translation copy of $A + B$ where it serves as a vertex of the convex hull $F$. By virtue of the ``boiling down'' construction of $R_\text{Win}$, all points of $\disc F \setminus \{\mathbf{v}\}$ must have been assigned values already. Call this assignment $\eta_\text{Step}$.

Observe that, in the proof of Lemma \ref{zs}, the initial point to be taken away is chosen arbitrarily among the vertices of $\lhull A \rhull + \lhull B \rhull$. Hence, there exists a smallest percolating seed $S_\text{Sum}$ for $F$ with $\mathbf{v} \not \in S_\text{Sum}$. Since $\Gamma$ is consistent over $\lhull A \rhull + \lhull B \rhull$, there also exists a solution $\eta_\text{Sum}$ to $\Gamma$ over $F$ where the points of $S_\text{Sum}$ are assigned the new formal indeterminates $y_1$, $y_2$, $\ldots$, $y_s$. Finally, let $f_1$, $f_2$, $\ldots$, $f_s$ be the rational functions drawn out of $\mathbb{F}(x_1, x_2, \ldots, x_s)$ which $\eta_\text{Step}$ assigns to the points of $S_\text{Sum}$.

Since $\Gamma$ is proper over $\mathbb{Z}^2$, by chaining together two suitable computational runs we find that substituting each $y_i$ with $f_i$ in the values of $\eta_\text{Sum}$ yields well-defined rational functions. These rational functions form a solution $\eta^\star_\text{Step}$ to $\Gamma$ over $F$. Furthermore, $\eta_\text{Step}$ and $\eta^\star_\text{Step}$ must agree over $\disc F \setminus \{\mathbf{v}\}$. Thus we may safely assign to $\mathbf{v}$ the value $\eta^\star_\text{Step}(\mathbf{v})$. \end{myproof}

One might rightfully wonder if Theorem \ref{cp} applies in any natural settings. We show below that it does establish consistency for one reasonably wide class of systems $\Gamma$. Other examples will be given in Section \ref{sd}.

\begin{proposition} \label{cpl} In the setting of Theorem \ref{cp}, suppose that the specifying polynomials of $\Gamma$ are linear and homogeneous. Then $\Gamma$ is indeed consistent on the plane. \end{proposition} % consistency in the plane for linear systems

\begin{myproof} It suffices to show that $\Gamma$ satisfies an analogue of Lemma \ref{lc}. We proceed to adapt its proof. Let $P = \sum_{\mathbf{a} \in A} a_\mathbf{a}x_\mathbf{a}$ and $Q = \sum_{\mathbf{b} \in B} b_\mathbf{b}x_\mathbf{b}$ be the specifying polynomials of $\Gamma$, and define the collections of polynomials $\mathcal{P}$ and $\mathcal{Q}$ as in the proof of Lemma \ref{lc}. We do not sum the members of $\mathcal{P}$ directly this time around, but instead for each $\mathbf{b} \in B$ we weight the image of $P$ upon translation by $\mathbf{b}$ with the coefficient $b_\mathbf{b}$ of $x_\mathbf{b}$ in $Q$; and similarly for $\mathcal{Q}$. The desired result now follows because both weighted sums work out to $\sum a_\mathbf{a}b_\mathbf{b}x_{\mathbf{a} + \mathbf{b}}$ over all $\mathbf{a} \in A$ and $\mathbf{b} \in B$. \end{myproof}

%% file: bpde-06-sd.tex
\section{Squares and Diamonds} \label{sd}

Here, we study some nonlinear examples of systems of two partial difference equations which are consistent over $\mathbb{Z}^2$.

Let $M^\square_k$ be the matrix of size $k \times k$, indexed by $1 \le i \le k$ and $1 \le j \le k$, with entry $x_{(i, j)}$ at position $(i, j)$. Let also $M^\lozenge_k$ be the matrix of size $k \times k$, indexed by $1 \le i \le k$ and $1 \le j \le k$, with entry $x_{(i - j, i + j)}$ at position $(i, j)$. Consider the polynomials $\det M^\square_k$ and $\det M^\lozenge_k$. Their patterns are a square grid and a diamond grid, respectively.

We write $\Phi^\square_k$ for the partial difference equation specified by $\det M^\square_k$, and $\Phi^\lozenge_k$ for the partial difference equation specified by $\det M^\lozenge_k$. We also write $\sd(m, n)$ for the system of two partial difference equations specified by the collection of polynomials $\{\det M^\square_m, \det M^\lozenge_n\}$. Since all of our specifying polynomials are with integer coefficients, in our solutions to these partial difference equations and systems of them we can draw the values out of an arbitrary field $\mathbb{F}$.

\begin{proposition} \label{sdc} Each one of $\sd(2, 2)$, $\sd(2, 3)$, $\sd(3, 2)$ is proper and consistent over~$\mathbb{Z}^2$. \end{proposition} % squares and diamonds consistency

\begin{myproof} Any nonzero constant assignment is a proper solution to $\sd(2, 2)$. Both of $\sd(2, 3)$ and $\sd(3, 2)$ admit proper solutions on the plane with values drawn out of the finite field $\mathbb{F}_5$: For $\sd(2, 3)$, to the point $(i, j)$ we assign $3$ when $i \equiv 0 \pmod 3$ and $1$ otherwise. For $\sd(3, 2)$, to the point $(i, j)$ we assign $3$ when $i \equiv j \pmod 3$ and $1$ otherwise.

We may now apply Theorem \ref{cp}. The consistency of each system over the Minkowski sum of the convex hulls of the patterns of its specifying polynomials can be verified by direct computation, as outlined near the beginning of Section \ref{cons-ii}. \end{myproof}

(For $\sd(2, 2)$, we can more generally say the following: Suppose, in the setting of Theorem~\ref{cp}, that each specifying polynomial of $\Gamma$ is of the form $\prod_{\mathbf{u} \in U} x_\mathbf{u} - \prod_{\mathbf{v} \in V} x_\mathbf{v}$, with $U$ and $V$ being two disjoint finite subsets of $\mathbb{Z}^2$. Then $\Gamma$ is indeed consistent on the plane. This is an ``exponentiated'' corollary of Proposition \ref{cpl}, as becomes apparent when we rewrite each one of the specifying polynomials of $\Gamma$ in the form $\prod_{\mathbf{u} \in U} x_\mathbf{u} \cdot \prod_{\mathbf{v} \in V} x_\mathbf{v}^{-1} = 1$. The proof is analogous to that of Proposition \ref{cpl}, with addition replaced by multiplication.)

On the other hand, the systems $\sd(m, n)$ with $2 \le m \le 10$ and $2 \le n \le 10$ are inconsistent over $\mathbb{Z}^2$ in all cases except for the ones of Proposition \ref{sdc}. Is this the full list of all positive integers $m \ge 2$ and $n \ge 2$ such that the system $\sd(m, n)$ is proper and consistent over $\mathbb{Z}^2$?

The systems $\sd(m, n)$ with $m = 2$ or $n = 2$ arise naturally in connection with Somos sequences, via the notion of diamond rank. We go on now to elaborate. The sequence $u$ in $\mathbb{F}$, indexed by $\mathbb{Z}$ and with nonzero terms, is said to satisfy a \emph{Somos recurrence} of order $k$ when there exist constants $a_0$, $a_1$, $\ldots$, $a_{\lfloor k/2 \rfloor}$ in $\mathbb{F}$, with $a_0 \neq 0$, such that for all $i$ the following condition holds: \[a_0u_iu_{i + k} + a_1u_{i + 1}u_{i + k - 1} + \cdots + a_{\lfloor k/2 \rfloor}u_{i + \lfloor k/2 \rfloor}u_{i + \lceil k/2 \rceil} = 0.\]

Of course, a Somos sequence $u$ of order $k$ is uniquely specified by its \emph{initial values} $u_1$, $u_2$, $\ldots$, $u_k$ and its \emph{coefficients} $a_0$, $a_1$, $\ldots$, $a_{\lfloor k/2 \rfloor}$. We define the \emph{free} Somos sequence of order $k$ to be the one with initial values $x_1$, $x_2$, $\ldots$, $x_k$ and coefficients $1$, $\alpha_1$, $\alpha_2$, $\ldots$, $\alpha_{\lfloor k/2 \rfloor}$, the $x$'s and $\alpha$'s being new formal indeterminates. Clearly, all of its terms will be integer-coefficient rational functions.

Given a sequence $u$ in $\mathbb{F}$, we write $u_\text{Even}$ and $u_\text{Odd}$ for the two ``halves'' of $u$ with $u_\text{Even}(i) = u(2i)$ and $u_\text{Odd}(i) = u(2i + 1)$. Given two sequences $u$ and $v$ in $\mathbb{F}$, each indexed by $\mathbb{Z}$, we write $u \times v$ for the arrangement over $\mathbb{Z}^2$ which assigns the value $u_iv_j$ to the point $(i, j)$.

Given an arrangement $\eta$ over $\mathbb{Z}^2$ with values drawn out of $\mathbb{F}$, we write $\operatorname{Dia}_\text{Even}(\eta)$ for the arrangement over $\mathbb{Z}^2$ which assigns the value $\eta(i - j, i + j)$ to the point $(i, j)$. Similarly, we write $\operatorname{Dia}_\text{Odd}(\eta)$ for the arrangement over $\mathbb{Z}^2$ which assigns the value $\eta(i - j, i + j + 1)$ to the point~$(i, j)$.

Notice that we may view $\eta$ as a matrix, with its rows and its columns indexed by $\mathbb{Z}$. By the \emph{rank} of $\eta$, we mean the rank of this matrix over $\mathbb{F}$. By the \emph{diamond rank} of $\eta$, we mean the maximum of the ordinary ranks of $\operatorname{Dia}_\text{Even}(\eta)$ and $\operatorname{Dia}_\text{Odd}(\eta)$.

It is well-known that, if $u$ is an order-$4$ Somos sequence, then $u \times u$ is of diamond rank at most~$2$. Similarly, if $u$ is an order-$5$ Somos sequence, then $u_\text{Even} \times u_\text{Odd}$ is of diamond rank at most~$2$. Furthermore, the terms of the free Somos sequences of orders $4$ and $5$ are Laurent polynomials. We point readers towards \cite{Be} for an in-depth discussion of this topic, including a bibliography. Meanwhile, we proceed to show how two generalisations of these results -- namely, Propositions \ref{mr} and \ref{mi} below -- can be derived quite cheaply from Proposition~\ref{sdc}.

We begin with one simple linear-algebraic observation, a variant of which appears also in \cite{Be}. The proof given there is rather more direct. However, we prefer here to demonstrate our newly developed toolbox in action.

\begin{lemma} \label{rank} Let $\eta$ be any proper solution to $\Phi^\square_{k + 1}$ over $\mathbb{Z}^2$. Then $\eta$ is of rank $k$. \end{lemma}

\begin{myproof} Fix any basis $\{\mathbf{u}_1, \mathbf{u}_2, \ldots, \mathbf{u}_k\}$ of $\mathbb{F}^k$. For each $j \in \mathbb{Z}$, there exists a unique $\mathbf{v}_j \in \mathbb{F}^k$ with $(\eta(1, j), \eta(2, j), \ldots, \eta(k, j)) = (\mathbf{u}_1 \cdot \mathbf{v}_j, \mathbf{u}_2 \cdot \mathbf{v}_j, \ldots, \mathbf{u}_k \cdot \mathbf{v}_j)$.

Observe that $\det M^\square_k$ is the denominator at $x_{(k + 1, k + 1)}$ of the polynomial $\det M^\square_{k + 1}$. Since $\eta$ is proper, we get that $\det M^\square_k(\eta) \neq 0$, and so $\{\mathbf{v}_1, \mathbf{v}_2, \ldots, \mathbf{v}_k\}$ is also a basis of $\mathbb{F}$. Hence, for each $i \in \mathbb{Z}$, there exists a unique $\mathbf{u}_i$ with $(\eta(i, 1), \eta(i, 2), \ldots, \eta(i, k)) = (\mathbf{u}_i \cdot \mathbf{v}_1, \mathbf{u}_i \cdot \mathbf{v}_2, \ldots, \mathbf{u}_i \cdot \mathbf{v}_k)$.

Consider now the union $S_\boxplus$ of the two infinite strips of width $k$ given by $\{(i, j) \mid i \in \mathbb{Z} \text{ and } 1 \le j \le k\}$ and $\{(i, j) \mid 1 \le i \le k \text{ and } j \in \mathbb{Z}\}$. It is a percolating seed for $\mathbb{Z}^2$ in the context of $\Phi^\square_{k + 1}$. Consider next the arrangement $\eta^\star$ over $\mathbb{Z}^2$ where $\eta^\star(i, j)$ is given by the scalar product $\mathbf{u}_i \cdot \mathbf{v}_j$. Since $\eta^\star$ is of rank $k$, it constitutes a solution to $\Phi^\square_{k + 1}$ over $\mathbb{Z}^2$.

Since $S_\boxplus$ percolates and $\eta$ is proper, we get that $\eta$ is the unique solution to $\Phi^\square_{k + 1}$ over $\mathbb{Z}^2$ which coincides with $\eta$ over $S_\boxplus$. However, $\eta$ and $\eta^\star$ certainly do coincide over $S_\boxplus$. Thus $\eta$ and $\eta^\star$ must also coincide globally, over $\mathbb{Z}^2$. \end{myproof}

By Lemma \ref{rank}, similarly each proper solution to $\Phi^\lozenge_{k + 1}$ over $\mathbb{Z}^2$ is of diamond rank $k$.

We see now that there exists a kind of reciprocity between the systems $\sd(m, n)$ and $\sd(n, m)$. Concretely, if $\eta$ is a proper solution to $\sd(m, n)$ over $\mathbb{Z}^2$, then both of $\operatorname{Dia}_\text{Even}(\eta)$ and $\operatorname{Dia}_\text{Odd}(\eta)$ must be solutions to $\sd(n, m)$ over $\mathbb{Z}^2$.

Consider next any proper solution $\eta$ to $\sd(2, k + 1)$ over $\mathbb{Z}^2$. By Lemma \ref{rank} for the square component of $\sd(2, k + 1)$, we get that $\eta$ is of unit rank. Hence, there exist two sequences $u$ and $v$ in $\mathbb{F}$, each indexed by $\mathbb{Z}$, with $\eta = u \times v$. On the other hand, by Lemma \ref{rank} for the diamond component of $\sd(2, k + 1)$, we get that $u$ and $v$ must satisfy a great wealth of determinant identities. This is our motivation for the following definition:

We say that two sequences $u$ and $v$ in $\mathbb{F}$, both indexed by $\mathbb{Z}$ and with nonzero terms, satisfy a \emph{mutual Somos recurrence} of order $k$ when there exist constants $a_0$, $a_1$, $\ldots$, $a_k$ and $b_0$, $b_1$, $\ldots$, $b_k$ in $\mathbb{F}$, with $a_0 \neq 0$ and $b_0 \neq 0$, such that for all $i$ the following conditions hold: \begin{align*} a_0u_iv_{i + k} + a_1u_{i + 1}v_{i + k - 1} + \cdots + a_ku_{i + k}v_i &= 0\\ b_0u_{i + k + 1}v_i + b_1u_{i + k}v_{i + 1} + \cdots + b_ku_{i + 1}v_{i + k} &= 0. \end{align*}

Of course, two mutually Somos sequences $u$ and $v$ of order $k$ are uniquely specified by their \emph{initial values} $u_1$, $u_2$, $\ldots$, $u_{k + 1}$, $v_1$, $v_2$, $\ldots$, $v_k$ and their \emph{coefficients} $a_0$, $a_1$, $\ldots$, $a_k$, $b_0$, $b_1$, $\ldots$, $b_k$. We define the two \emph{free} mutually Somos sequences of order $k$ to be the ones with initial values $x_1$, $x_2$, $\ldots$, $x_{k + 1}$, $y_1$, $y_2$, $\ldots$, $y_k$ and coefficients $1$, $\alpha_1$, $\alpha_2$, $\ldots$, $\alpha_k$, $1$, $\beta_1$, $\beta_2$, $\ldots$, $\beta_k$, the $x$'s, $y$'s, $\alpha$'s, and $\beta$'s being new formal indeterminates. Clearly, all of their terms will be integer-coefficient rational functions.

The ordinary Somos recurrence is a special case of the mutual Somos recurrence. Consider first the setting where $b_0 = 1$, $b_1 = -1$, $b_2 = b_3 = \cdots = b_k = 0$, $v_1 = u_{k + 1}$. Then $v_i = u_{i + k}$ for all $i$, and so $u$ becomes the Somos sequence of order $2k$ with coefficients $a_0$, $a_1$, $\ldots$, $a_k$ and initial values $u_1$, $u_2$, $\ldots$, $u_{k + 1}$, $v_2$, $v_3$, $\ldots$, $v_k$. Consider next the setting where $a_i = b_i$ for all~$i$. Then the interleaved sequence $\ldots$, $u_1$, $v_1$, $u_2$, $v_2$, $u_3$, $v_3$, $\ldots$ becomes the Somos sequence of order $2k + 1$ with coefficients $a_0$, $a_k$, $a_1$, $a_{k - 1}$, $\ldots$, $a_{\lceil k/2 \rceil}$ and initial values $u_1$, $v_1$, $u_2$, $v_2$, $\ldots$, $v_k$,~$u_{k + 1}$. Or, in summary, the mutual Somos recurrence of order $k$ simulates the ordinary Somos recurrences of orders both $2k$ and $2k + 1$.

We get that, in particular, the mutual Somos recurrence of order $2$ simulates the ordinary Somos recurrences of orders both $4$ and $5$. This is how the aforementioned properties of the order-$4$ and order-$5$ Somos recurrences become corollaries of Propositions \ref{mr}~and~\ref{mi} below.

We return now to $\sd(2, k + 1)$ and its proper solution $\eta = u \times v$ over $\mathbb{Z}^2$. It is straightforward to see that $u$ and $v$ must satisfy a mutual Somos recurrence of order $k$. For example, if $k = 2$, the coefficients of this mutual Somos recurrence will be given by \begin{align*} a_0 &= \det \begin{pmatrix} \eta(2, 2) & \eta(3, 1)\\ \eta(3, 3) & \eta(4, 2) \end{pmatrix} & a_1 &= \det \begin{pmatrix} \eta(3, 1) & \eta(1, 3)\\ \eta(4, 2) & \eta(2, 4) \end{pmatrix} & a_2 &= \det \begin{pmatrix} \eta(1, 3) & \eta(2, 2)\\ \eta(2, 4) & \eta(3, 3) \end{pmatrix}\\ b_0 &= \det \begin{pmatrix} \eta(2, 3) & \eta(3, 2)\\ \eta(3, 4) & \eta(4, 3) \end{pmatrix} & b_1 &= \det \begin{pmatrix} \eta(4, 1) & \eta(2, 3)\\ \eta(5, 2) & \eta(3, 4) \end{pmatrix} & b_2 &= \det \begin{pmatrix} \eta(3, 2) & \eta(4, 1)\\ \eta(4, 3) & \eta(5, 2) \end{pmatrix}. \end{align*}

Since the system $\sd(2, 3)$ is consistent over $\mathbb{Z}^2$ by Proposition \ref{sdc}, for it the latter observation admits a sort of converse:

\begin{proposition} \label{mr} Let $u$ and $v$ be two mutually Somos sequences of order $2$. Then $u \times v$ is of diamond rank at most $2$. So, in particular, $u \times v$ is a solution to $\sd(2, 3)$ over $\mathbb{Z}^2$. \end{proposition} % mutual rank

\begin{myproof} Observe that $S_\boxtimes = \{(1, j) \mid 1 \le j \le 4\} \cup \{(i, 1) \mid 2 \le i \le 5\}$ is a smallest percolating seed for $\mathbb{Z}^2$ in the context of $\sd(2, 3)$. Let $\eta_\boxtimes$ be the solution to $\sd(2, 3)$ over $\mathbb{Z}^2$ where the points of $S_\boxtimes$ are assigned the new formal indeterminates $z_1$, $z_2$, $\ldots$, $z_8$.

It suffices to consider the case of the two free mutually Somos sequences $u_\text{Free}$ and $v_\text{Free}$ of order $2$. Let $\sigma$ be the substitution which, for each point $(i, j)$ of $S_\boxtimes$, replaces its associated $z_k$ with the rational function $u_\text{Free}(i)v_\text{Free}(j)$. Since setting all indeterminates in $u_\text{Free}$ and $v_\text{Free}$ to unity yields two mutually Somos sequences in $\mathbb{F}_5$ whose product is a proper solution to $\sd(2, 3)$ over $\mathbb{Z}^2$, we get that $\sigma$ maps each value of $\eta_\boxtimes$ onto a well-defined rational function. Furthermore, these rational functions form a solution $\eta^\star_\boxtimes$ to $\sd(2, 3)$ over $\mathbb{Z}^2$.

Clearly, $\eta^\star_\boxtimes$ factors into two mutually Somos sequences whose initial values coincide with those of $u_\text{Free}$ and $v_\text{Free}$. By a straightforward calculation, their coefficients also coincide with the coefficients of $u_\text{Free}$ and $v_\text{Free}$. We conclude that, in fact, $\eta^\star_\boxtimes = u_\text{Free} \times v_\text{Free}$. \end{myproof}

One useful lemma in \cite{Be} shows that, in the setting where $\mathbb{F}$ is the fraction field of some unique factorisation domain $\mathbb{U}$, the irreducibles which can occur in the denominators of the terms of $u$ are narrowly constrained when $u \times u$ is of finite diamond rank. The result is not too difficult to generalise with the product $u \times v$ of two distinct sequences $u$ and $v$. By means of this generalisation, Proposition \ref{mr} yields the following corollary:

\begin{proposition} \label{mi} The terms of the two free mutually Somos sequences of order $2$ are Laurent polynomials in $\mathbb{Z}[\alpha_1, \alpha_2, \beta_1, \beta_2, x_1, x_1^{-1}, x_2, x_2^{-1}, x_3, x_3^{-1}, y_1, y_1^{-1}, y_2, y_2^{-1}]$. \end{proposition} % mutual integrality

One might naturally wonder if anything similar is going on with the Somos recurrences of higher orders. For example, it is well-known that if $u$ is an order-$6$ Somos sequence, then $u \times u$ is of diamond rank at most $4$. Similarly, if $u$ is an order-$7$ Somos sequence, then $u_\text{Even} \times u_\text{Odd}$ is of diamond rank at most $4$. (Once again, details and references can be found in \cite{Be}.) Do these properties admit a generalisation with the mutual Somos recurrence of order $3$? We make the following conjecture, which does capture the ordinary Somos recurrences of orders $6$ and $7$ as special cases:

\begin{conjecture} \label{mho} Let $u$ and $v$ be two mutually Somos sequences of order $3$ whose coefficients satisfy $a_0b_1b_2a_3 = b_0a_1a_2b_3$. Then $u \times v$ is of diamond rank at most $4$. So, in particular, $u \times v$ is a solution to $\sd(2, 5)$ over $\mathbb{Z}^2$. \end{conjecture} % mutual higher orders

Since none of the systems $\sd(m, n)$ with $\{m, n\} = \{2, 4\}$ or $\{2, 5\}$ are consistent over $\mathbb{Z}^2$, the method we employed in the proof of Proposition \ref{mr} is not applicable directly. On the other hand, the twinning method of \cite{Be} seems more promising.

The gist of it is as follows: First, we find suitable ``invariants'' for the mutual Somos recurrence of order $3$. Then, based on them, we introduce a notion of ``twinning'' for pairs of mutually Somos sequences of order $3$. Finally, we show that whenever two pairs of mutually Somos sequences $(u', v')$ and $(u'', v'')$ of order $3$ are twins, all of the products $u' \times u''$, $u' \times v''$, $v' \times u''$, $v' \times v''$ are of diamond rank at most $4$.

This strategy, if successful, would yield some noteworthy corollaries. Concretely, since each pair of mutually Somos sequences $(u, v)$ of order $3$ must necessarily be twinned with itself, it would follow that the products $u \times u$ and $v \times v$ are also of diamond rank at most $4$. So, in particular, each one of $u$ and $v$ would need to satisfy some ordinary Somos recurrence of even order at most $8$ as well as some ordinary Somos recurrence of odd order at most $9$.

%% file: bpde-07-further.tex
\section{Higher Dimensions and Further Work} \label{further}

Let $\mathcal{A}$ be a collection of patterns in $\mathbb{Z}^d$, as in the introduction. When does $\mathcal{A}$ admit a finite percolating seed? For example, is it true that a finite percolating seed exists if and only if the patterns of $\mathcal{A}$ do not share a stable direction?

Then, in the setting where a finite percolating seed does exist: What is the smallest size of such a seed, in terms of $\mathcal{A}$? It is unlikely that a single neat formula would apply universally. However, Theorems \ref{1d} and \ref{2d} hint at the potential relevance of the mixed volume with $d$ patterns in $d$ dimensions. We offer the following lower bound:

\begin{theorem} \label{dd} For $d$ patterns in $d$ dimensions, the size of the smallest percolating seed is bounded from below by the mixed volume of the patterns' convex hulls. \end{theorem}

We already know that, in one and two dimensions, the lower bound of Theorem \ref{dd} is attained whenever $\mathcal{A}$ admits a finite percolating seed. Do analogous results hold in any higher dimensions as well?

While our proofs of Theorems \ref{1d} and \ref{2d} were self-contained and purely combinatorial, for Theorem \ref{dd} we must instead invoke a bit of advanced algebraic-geometric machinery. Let $\mathcal{A} = \{A_1, A_2, \ldots, A_d\}$. We write $\mathbf{x}^\mathbf{a}$ for $x_1^{a_1}x_2^{a_2} \cdots x_d^{a_d}$, where $x_1$, $x_2$, $\ldots$, $x_d$ are new formal indeterminates and $\mathbf{a} = (a_1, a_2, \ldots, a_d) \in \mathbb{Z}^d$. Consider the collection of Laurent polynomials $\mathcal{E} = \{\sum_{\mathbf{a} \in A} c_{\mathbf{a}, A}\mathbf{x}^\mathbf{a} \mid A \in \mathcal{A}\}$, with complex coefficients $c_{\mathbf{a}, A}$.

Bernstein's theorem \cite{CLOS} states that, for generic $c_{\mathbf{a}, A}$, the number of roots to $\mathcal{E}$ with nonzero complex components equals $\mix(\lhull A_1 \rhull, \lhull A_2 \rhull, \ldots, \lhull A_d \rhull)$. Here, the qualifier ``generic'' means that there exists a complex-coefficient polynomial $Q$, with $|A_1| + |A_2| + \cdots + |A_d|$ arguments and depending solely on $\mathcal{A}$, such that the conclusion holds for all complex coefficients $c_{\mathbf{a}, A}$ which do not cause $Q$ to vanish.

\begin{myproof}[Proof of Theorem \ref{dd}] Fix any nonzero complex coefficients $c_{\mathbf{a}, A}$ generic in the sense just described, and consider the system of partial difference equations $\Gamma$ specified by the collection of linear and homogeneous polynomials $\{\sum_{\mathbf{a} \in A} c_{\mathbf{a}, A}x_\mathbf{a} \mid A \in \mathcal{A}\}$.

Let $\mathbf{e} = (e_1, e_2, \ldots, e_d) \in (\mathbb{C} \setminus \{0\})^d$ be any root to $\mathcal{E}$ with nonzero complex components. Then the arrangement $\eta$ given by $\eta(\mathbf{a}) = \mathbf{e}^\mathbf{a}$ is a solution to $\Gamma$. Furthermore, all such solutions to $\Gamma$ are linearly independent. Indeed, it is a standard result in multivariate interpolation that with pairwise distinct $\mathbf{e}_1$, $\mathbf{e}_2$, $\ldots$, $\mathbf{e}_k$, the associated arrangements $\eta_1$, $\eta_2$, $\ldots$, $\eta_k$ are already linearly independent when restricted to the subregion $\{(a_1, a_2, \ldots, a_d) \mid a_i \ge 0 \text{ for all } i \text{ and } a_1 + a_2 + \cdots + a_d < k\}$ of $\mathbb{Z}^d$.

We conclude that $\dim \sol \Gamma$ is bounded from below by $\mix(\lhull A_1 \rhull, \lhull A_2 \rhull, \ldots, \lhull A_d \rhull)$. On the other hand, by a straightforward generalisation of Lemma \ref{dbps}, the former quantity itself bounds from below the smallest size of a percolating seed for $\mathcal{A}$. \end{myproof}

How might we approach a proof that the lower bound of Theorem \ref{dd} is in fact attained, for some class of pattern collections $\mathcal{A}$? The statements of Lemmas \ref{zs} and \ref{se} in the proof of Theorem~\ref{2d} suggest the following two-step strategy: First, we show that there exists a percolating seed of size $\mix(\lhull A_1 \rhull, \lhull A_2 \rhull, \ldots, \lhull A_d \rhull)$ for the region $\lhull A_1 \rhull + \lhull A_2 \rhull + \cdots + \lhull A_d \rhull$. Second, we show that the integer points of the latter region themselves form a percolating seed for the full space $\mathbb{Z}^d$. Still, it is far from obvious how the supporting discrete-geometric machinery of Section \ref{plane} could be adapted to higher dimensions.

We move on now to the interplay between bootstrap percolation and partial difference equations. The toolbox developed in Sections \ref{cons-i} and \ref{cons-ii} can be generalised so as to handle many non-multilinear $\Gamma$. We proceed next to sketch this generalisation in broad strokes.

Let $\mathcal{A} = \{A_1, A_2, \ldots, A_k\}$. We introduce a second collection of patterns $\mathcal{A}^\star = \{A^\star_1, A^\star_2,\allowbreak \ldots, A^\star_k\}$, too, with $A^\star_i \subseteq A_i$ for all $i$. We say that $\mathcal{A}^\star$ \emph{marks} $\mathcal{A}$, with the marked points of $A_i$ being $A^\star_i$.

Let $\tau$ be any translation. Suppose that the point $\mathbf{a}$ completes the translation copy $\tau(A_i)$ of some member of $\mathcal{A}$. We say that the completion is \emph{legal} when $\mathbf{a} \in \tau(A^\star_i)$; i.e., when $\mathbf{a}$ plays the role of a marked point within $\tau(A_i)$.

We define the marked bootstrap percolation game specified by $\mathcal{A}$ and $\mathcal{A}^\star$ similarly to the ordinary bootstrap percolation game specified by $\mathcal{A}$, with one crucial modification: We are only allowed to occupy a new point when it completes a translation copy of some member of $\mathcal{A}$ in a legal manner. Notice that this setup is equivalent to the framework of ``$\mathcal{U}$-bootstrap percolation'' as set forth in \cite{BSU}.

Let $\Gamma$ be the system of partial difference equations specified by the collection of polynomials $\mathcal{P} = \{P_1, P_2, \ldots, P_k\}$, with $A_i$ being the pattern of $P_i$ for all $i$. We now require that each $P_i$ must be linear merely in its marked indeterminates; i.e., in those indeterminates $x_\mathbf{a}$ for which $\mathbf{a} \in A^\star_i$. This allows us to consider non-multilinear $\Gamma$, too, with the ordinary bootstrap percolation game specified by $\mathcal{A}$ replaced with the marked bootstrap percolation game specified by $\mathcal{A}$ and $\mathcal{A}^\star$.

Notice that, in this setting, propriety becomes a property not just of the system and the region, but also of our chosen markings. We are not forced to mark every single point $\mathbf{a}$ in $A_i$ such that $P_i$ is linear in $x_\mathbf{a}$, but can instead select among them.

Furthermore, we cannot afford any longer to define consistency solely in terms of path-independence, as some newly occupied points might bring about illegal completions. Instead, we must define $\Gamma$ to be consistent from $S$ to $T$ when, for every run $R$ from $S$ to $T$, the free computation associated with $R$ assigns identical values to identical points; and, besides, these values always form a solution to $\Gamma$ over $T$.

Other than that, the material of Sections \ref{cons-i} and \ref{cons-ii} carries over with no significant difficulties.

The naive variants of Theorems \ref{1d}, \ref{2d}, \ref{cp}, and \ref{ce} in the marked setting would state simply that all points must be marked. However, careful examination of the proofs reveals that in reality we can get away with much less. Concretely, all four of Theorems \ref{1d}, \ref{2d}, \ref{cp}, and \ref{ce} continue to hold in the marked setting when it is merely the vertices of the patterns' convex hulls that are marked. These are substantial generalisations, with the strengthened variants of Theorems \ref{cp} and \ref{ce} applying to vast swathes of non-multilinear $\Gamma$.

What would be some other exciting examples of systems of two partial difference equations in two dimensions consistent over $\mathbb{Z}^2$, apart from the ones identified in Sections \ref{cons-ii} and \ref{sd}? One potential approach might be to run exhaustive searches over all such systems where the sizes and bounding boxes of the patterns as well as the degrees and monomial counts of the specifying polynomials are bounded from above somehow.

Finally, it would be interesting to see some higher-dimensional criteria for consistency similar to Theorem \ref{cp}. The preceding discussion indicates that the setting of $d$ partial difference equations in $d$ dimensions, with the patterns of the specifying polynomials admitting a finite percolating seed, could be especially fertile ground for further investigations.